\documentclass [12pt,a4paper]{article}
\usepackage{epsfig,amsmath,amssymb}
\usepackage{natbib}
\usepackage{bm}
\usepackage{psfrag}
\usepackage{wrapfig}
\usepackage{booktabs}
\usepackage{geometry}
\usepackage{pdflscape}
\def\e{\hbox{E}}

\def\min{\hbox{min}}

\def\cdf{\hbox{cdf\;}}
\def\df{\hbox{df\;}}
\def\kld{\hbox{KLD\,}}
\def\pdQ{\hbox{pdQ\,}}
\newcommand{\D}{\mathcal{D}\,}
\newcommand{\F}{\mathcal{F}\,}
\newcommand{\U}{\mathcal{U}}
\newtheorem{definition}{Definition}
\newtheorem{proposition}{Proposition}

\begin{document}

\title{The Shapes of Things to Come: \\ Probability Density Quantiles}
\author{Robert G. Staudte}
\date{16 December, 2016}
\maketitle

\begin{abstract}
For every discrete or continuous location-scale family having a square-integrable density,
there is a unique continuous probability distribution on the unit interval that is determined by the density-quantile composition introduced by Parzen in 1979. These
{\em probability} density quantiles (\pdQ s) only differ in {\em shape}, and can be usefully
compared with  the Hellinger distance or  Kullback-Leibler divergences. Convergent empirical estimates of these \pdQ s are provided, which leads to a robust global fitting procedure of
shape families to data. Asymmetry can be measured in terms of distance or divergence of \pdQ s from the
symmetric class. Further, a precise classification of shapes by tail behavior can be  defined simply in terms of \pdQ boundary derivatives.
\end{abstract}

{\em Keywords:  \ asymmetry; Hellinger metric; quantile density; Kullback-Leibler divergence; Lagrange multipliers; tail-weight}

\section{Introduction}\label{sec:intro}
   The study of shapes of probability distributions is simplified by viewing them
through the revealing composition of density with quantile function. Following normalization,
the resulting functions are not only location-scale free, but densities of absolutely continuous
distributions having the same support. This microcosm of {\em probability density quantiles} carries essential information regarding shapes and allows for simpler classification by asymmetry and
tail weights. It also leads to an alternative and effective method for fitting shape families to data.

\subsection{Background and summary}

In the seminal work \cite{par-1979} proposed that traditional statistical inference be connected to exploratory data
analysis through transformations from standard continuous models (normal, exponential)
to the realm of density quantile functions. If the standard models were rejected by goodness-of-fit tests, the next step was nonparametric modeling, with the emphasis
on what could be gleaned from sample quantile functions and time series methods. The quantile approach to data analysis
was earlier championed by \cite{tukey-1962,tukey-1965,tukey-1977}, who also provided insightful commentary into Parzen's proposals.  Further work on quantile-based data modeling can be found in \cite{gilch-2000}, \cite{par-2004}, while \cite{jones-1992} investigates estimation of density quantile functions and their reciprocals.

Here we study the classification of shapes for {\em probability} density
quantiles.  While this class is limited by the requirement of square-integrability of the density function, it is rich enough to warrant investigation because the transformation from density to the normalized composition of density function with quantile function allows for comparison of shapes of both discrete and continuous models using the Hellinger metric and Kullback-Leibler divergences---with non-trivial and informative results.

In this Section we formally introduce \pdQ s for discrete and continuous distributions and provide numerous examples, including moments in the continuous case. Given data, in Section~\ref{sec:epdQ} we describe  methods for estimating the pdQ's of discrete and continuous distributions, respectively. These are employed in Section~\ref{sec:applics} where global robust
fitting of shape parameter families to data  based on the Hellinger distance from the empirical \pdQ are implemented.

In Section~\ref{sec:asymmetry} we measure asymmetry of the \pdQ \, in terms of its distance or divergence from the class of symmetric distributions; and, we show they can be predicted by the skewness coefficient of the \pdQ .\ In Section~\ref{sec:tailwt} we introduce a simple but absolute tail-weight classification in terms of the boundary derivatives of the \pdQ s.  Further challenges are posed in Section~\ref{sec:summary}.

\subsection{Definitions and properties}\label{defns}

 Let $\F $ be the class of all right-continuous cumulative distribution functions (\cdf s) on the real line.  For each $F\in \F$ define the associated left-continuous {\em quantile function} of $F$ by $ Q(u)\equiv \inf \{x:\ F(x)\ge u \}$, for $0< u< 1.$  When the random variable $X$ has \cdf $F$, write $X\sim F$. In particular,
  let $U\sim \U $ where $U$ has the uniform distribution on [0,1].
Let $\F '=\{F\in \F :\; f=F' \text { exists and is positive}\}$.
          For each $F\in \F '$ we follow \cite{par-1979} and define  the {\em quantile density} function $q(u)=Q'(u)=1/f(Q(u)),$  \cite{tukey-1965} also recognized its importance and called it the {\em sparsity index}.  Its reciprocal $fQ(u)\equiv f(Q(u))$ is called  the {\em density quantile} function.
In order to convert density quantiles into {\em probability} densities we need to compute
  $\kappa \equiv \e [fQ(U)] = \int f^2(x)\,dx .$

\begin{definition}\label{def1}
 For $F\in \F '$, assume $\kappa =\e [fQ(U)]$ is finite; that is, $f$ is square integrable.
 Then define the {\em probability  density quantile} or \pdQ of $F$ by  $f^*(u)=fQ(u)/\kappa $, $0<u<1$.  Let $\F '^{*}\subset \F ' $ denote the class of all such $F$.
\end{definition}

Not all densities are square integrable, and for such densities a \pdQ does not exist.
Examples are the Chi-squared densities with degrees of freedom $\nu \leq 1$. Others are the Beta$(a,b )$ densities with $a\leq 1/2$ or $b\leq 1/2.$
Unless otherwise noted we follow standard definitions for distributions as described in \cite{J-K-B-1994,J-K-B-1995}.

  An important property of $\pdQ $s is that they are location-scale invariant.  For if  $F_{a,b}(\cdot )\equiv F((\cdot -a)/b)$ for arbitrary $a$ and $b>0$ defines the location-scale family generated by $F=F_{0,1}\in \F $, then the quantile function for $F_{a,b}$ is $Q_{a,b}(u)=a+b\,Q(u)$. Further, if $F\in \F '$ the quantile density is $q_{a,b}(u)=b\,q(u)$; thus the quantile density is location-invariant and scale equivariant.
 Clearly $f_{a,b}^*$ is also scale invariant, and one can write $f_{a,b}^*=f_{0,1}^*=f^*.$ Thus when comparing the graphs of different $f^*$s, we are comparing only their {\em shapes}.

Conversely, given an $f^*\in \F '^{*}$, one can identify the family $\{F_{a,b}:\ a,b>0\}$. For if it is known that $f^* $ has the form $f^* =(fQ)/\kappa $, for some unknown $F$ with associated density $f=F'$, inverse $Q$, quantile density $q=1/(fQ)$ and $\kappa =\e [fQ(U)]$, then one can
reconstruct $Q(u)=\int _0^u q(t)\,dt+c=\kappa \int _0^u \{f^*(t)\}^{-1}\,dt +c$; thus $Q$ is determined up to location and scale, as is $F$. An open question is whether, given an arbitrary continuous distribution with probability density $g$ on $(0,1),$ does there exist an $F\in \F '$ such that $g=f^*$?

Another property of \pdQ s is that they ignore flat spots in $F$. For example, the \pdQ \; $f^*_\text{gap}$
of  $f _\text{gap}(x)=(e^x/2) \;I_{(-\infty ,0]}(x) + I_{[1/2,1]}(x)$ equals that of $f(x)= (e^x/2) \;I_{(-\infty ,0]}(x) + I_{[0,1/2]}(x)$. Thus it is only the shape of the distribution on its {\em support} that is captured by the \pdQ .

\clearpage
\newpage

\begin{figure}[t!]
\begin{center}
\includegraphics[scale=0.9]{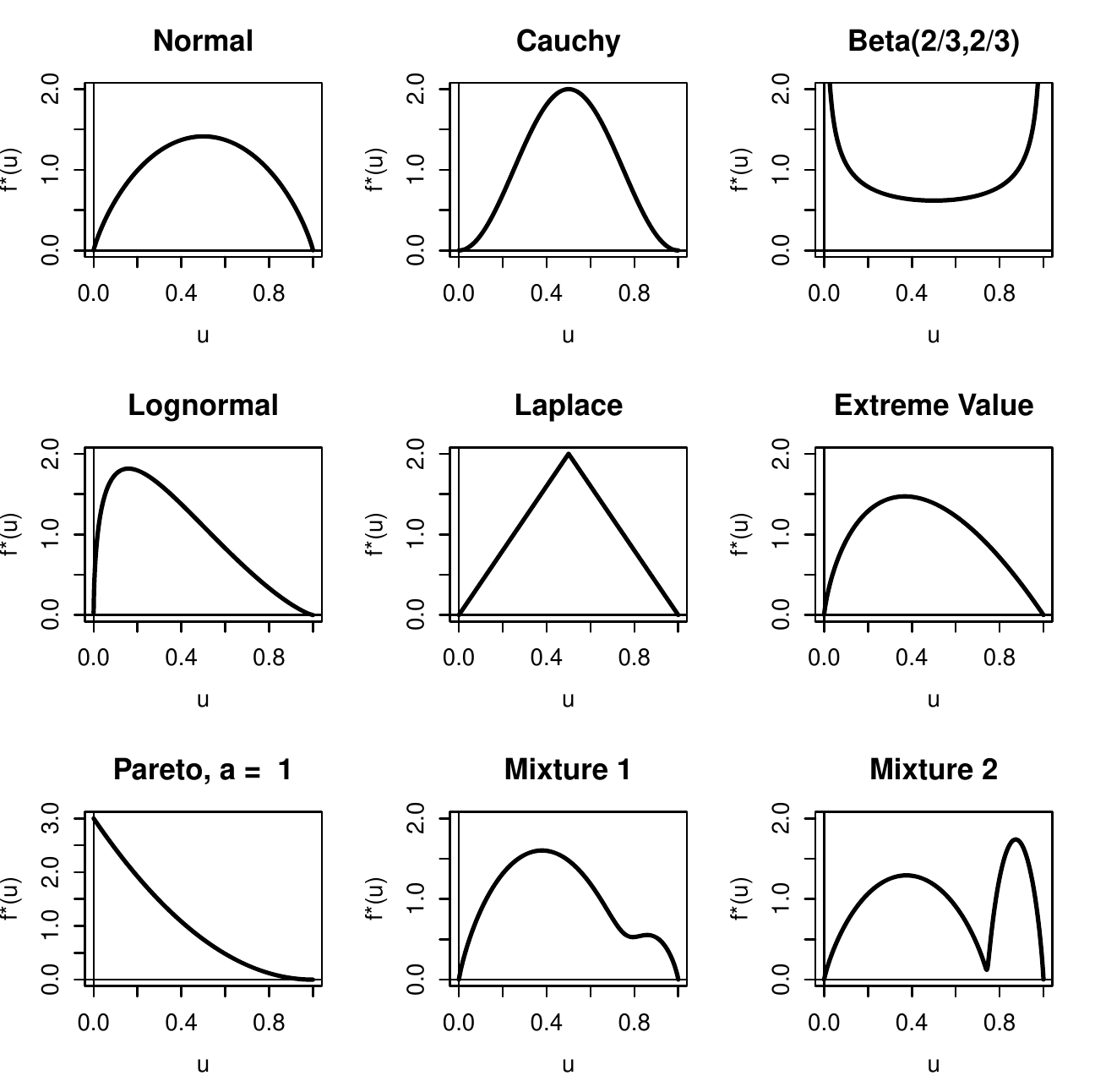}
\caption{\footnotesize  \em \pdQ s of some standard continuous distributions, plus two mixtures of
normal distributions.  Details are given in Section~\ref{sec:contpdQs}. \label{fig1}}
\end{center}
\end{figure}

\subsection{Examples of continuous \pdQ s}\label{sec:contpdQs}

 Many of the \pdQ s  $f^*(u)=fQ(u)/\kappa $ discussed in this section are merely normalized versions of density quantiles $fQ$ described in  \cite{par-1979}; some formulae for $f^*(u)$ are given in Table~\ref{table1}.

For the normal distribution $\Phi $ with density $\varphi $, we write $z_u=\Phi ^{-1}(u)$ and then it has  $\varphi ^*(u)=2\sqrt {\pi }\,\varphi (z_u).$ Its graph
is shown in the upper left plot of Figure~\ref{fig1}. It is quite close to that of a quadratic function which is symmetric about 1/2 and passes through $(0,0)$ and $(1/2, \sqrt 2)$, and not unlike the exact quadratic $f^*(u)=6u(1-u)$ corresponding to the logistic distribution (not shown).   The U-shaped Beta(2/3,2/3) distribution retains its U-shape, and the bell-shaped Cauchy  retains a bell-shape  after transformation to the quantile scale. The Laplace (double-exponential) family transforms to the symmetric triangular distribution.

The  Tukey($\lambda $) family is defined by its quantile function $Q_\lambda (u)=\{u^\lambda -(1-u)^\lambda \}/\lambda $ for $\lambda \neq 0$ and $Q_0(u)=\ln (u/(1-u))$, but
in general no closed form expression is available for its density. For $\lambda >0$ the density has finite
support $[-1/\lambda ,+1/\lambda ]$, and otherwise has infinite support. It has \pdQ \
$f^*_\lambda (u)=1/\{Q'_\lambda (u)\kappa _\lambda \}  =\{u^{\lambda -1} +(1-u)^{\lambda -1}\}^{-1}/\kappa _\lambda .$ The constant
$\kappa _\lambda $ required to make $\int f^*_\lambda (u)=1$ can be obtained by numerical integration but a good approximation is given by:
\begin{equation}\label{tukeylambdaK}
\kappa _\lambda ~\approx ~\left\{
       \begin{array}{ll}
      \frac {3^\lambda }{6}~, & \hbox{\quad } \qquad \lambda \leq 1~;\\
     \frac{\lambda}{2}     ~, & \hbox{ \quad } 1\leq  \lambda \leq 2    ~;  \\
      \left (\frac {2}{\pi }\right )^{\lambda -2}~, & \hbox{ \quad } 2 \leq  \lambda \leq 6~.
       \end{array}
     \right.
 \end{equation}
For $\lambda \leq 2$ the absolute error of this approximation is less than 0.005, and  for $2\leq \lambda \leq 6$ it is less 0.1. The absolute relative error of approximation is less than 0.06 for $-1\leq \lambda \leq 6.$

\clearpage
\newpage

The Pareto families (Type I or II) with shape parameter $a>0$ have \pdQ s defined by
$f^*_a (u)=(2+1/a)(1-u)^{1+1/a}$.
As $a\to \infty $ the graphs of $f^*_a$ rapidly approach a triangle, the graph of the exponential distribution \pdQ .
  The remaining graphs in Figure~\ref{fig1} are  3:1 mixtures of two normal distributions; Mixture 1 has components $N(0,1)$ and $N(3,1)$, and Mixture 2 components are $N(0,1)$ and $N(3,1/4^2)$.
In summary, all uniform densities transform to the standard uniform;  bell-shaped
densities with short tails transform roughly to quadratic functions; and exponential distributions correspond to triangular shapes. Mixtures of normal densities appear to transform to a  \lq mixture of quadratics\rq .\

\begin{table}[t!]
\begin{center}
\caption{{\bf Some continuous distributions $F$ and their quantile functions and \pdQ s.\ } \em In general, we denote $x_u=Q(u)=F^{-1}(u)$, but for the normal $F=\Phi $ with density $\varphi $, we write
$z_u=\Phi ^{-1}(u)$.  The power function distribution is also known as the Beta$(b,1)$ family. For the Weibull$(\beta )$ family, no simple approximation for $\kappa _\beta $ is available. Similarly for the Tukey$(\lambda )$ family, $\kappa _\lambda $ is not available, but a simple approximation is given in Section~\ref{sec:contpdQs}.  \label{table1}}
\begin{tabular}{lccccc}
                    &  $F(x)$       && $Q(u)$           &&  $f^*(u)$              \\
  \cline{1-6}\\[-.3cm]
Power$(b)$ &  $x^b,~0<x<1$        && $u^{1/b} $       &&   $(2-\frac{1}{b})u^{1-1/b}~,~ b>1/2$  \\[.2cm]
Laplace        &  $\qquad \ e^{x}/2,~ x< 0$ && $\,\quad \qquad \ln (2u),~u\leq 0.5$   &&  $2\,\min\{u,1-u\}  $   \\
                    &  $1-e^{-x}/2,~x\geq 0$ && $ -\ln (2(1-u)),~u\geq 0.5$   &&             \\[.2cm]
 Logistic           &  $e^x\,(1+e^x)^{-1}$      && $\ln(u/(1-u))$              &&  $6u(1-u)$            \\[.2cm]
 Extreme Value      &  $e^{-e^{-x}} $    && $-\ln(-\ln (u))$            &&  $-4u\ln(u)$          \\[.2cm]
 Cauchy             &  $\frac {1}{\pi}\arctan (x)-\frac {1}{2}$ &&  $\tan\{\pi(u-0.5)\}$  && $2\sin ^2(\pi u)$    \\[.2cm]
 Tukey$(\lambda )$   &  $-$  && $\lambda ^{-1}\{u^\lambda -(1-u)^\lambda \}$ && $
\{\kappa _{\lambda }\;(u^{\lambda -1} +(1-u)^{\lambda -1})\}^{-1}$         \\[.2cm]
 Normal             &  $\Phi (x)$         && $z_u$                       &&  $2\sqrt \pi\,\varphi(z_u)$ \\[.2cm]
Lognormal       &  $ \Phi (\ln (x)),~ x> 0 $  && $e^{z_u}$   &&  $\frac {2\sqrt \pi\;}{e^{1/4}}\, \varphi (z_u)\,e^{-z_u}$ \\[.2cm]
 Type I Pareto$(a)$        &  $1- x^{-a},~ x> 1 $     && $(1-u)^{-1/a}$         &&  $(2+\frac{1}{a})\,(1-u)^{1+1/a}$ \\[.2cm]
 Exponential        &  $1-e^{-x},~ x> 0 $           && $-\ln (1-u)$                &&  $2(1-u)$        \\[.2cm] Weibull$(\beta)$    &  $1- e^{-x^{\beta}},~ x> 0 $  &&$\{-\ln (1-u)\}^{1/\beta}$   &&  $\frac {\beta (1-u)}{\kappa _\beta }\,\{-\ln (1-u)\}^{1-1/\beta }$
 \end{tabular}
\end{center}
\end{table}

\begin{table}[h!]
\begin{center}
\caption{\label{table2}{\bf Examples of  skewness  and kurtosis for \pdQ s of continuous distributions $F$.}
 \em  For the symmetric cases on the left, $\mu ^*=0.5$ and $\gamma _1 ^*=0$.}
\begin{tabular}{ccccccccc}
$F$&$\sigma^*$ &$\gamma _2^*$ &$\quad $& $F$   & $\mu ^*$ & $\sigma^*$ & $\gamma _1 ^*$ & $\gamma _2^*$\\
\hline\\[-.1cm]
Beta(2/3,2/3)   &  0.3561 &  1.4846 &  &     Pareto(0.5) & 0.2000 & 0.1633 & 1.0498 & 3.6964 \\
Uniform         &  0.2887 &  1.8000 &  &     Pareto(1)   & 0.2500 & 0.1936 & 0.8607 & 3.0952 \\
Laplace         &  0.2041 &  2.4000 &  &     Pareto(2)   & 0.2857 & 0.2130 & 0.7318 & 2.7566 \\
Cauchy          &  0.1808 &  2.4062 &  &    Weibull(2)   & 0.4557 & 0.2393 & 0.1315 &  2.0714 \\
$t_2$           &  0.2041 &  2.2500 &  &    $\chi ^2_2$  & 0.3333 & 0.2357 & 0.5657 &  2.4000\\
$t_3$           &  0.2131 &  2.1961 &  &    $\chi ^2_3$  & 0.3849 & 0.2354 & 0.3808 &  2.2246\\
$t_5$           &  0.2207 &  2.1527 &  &    $\chi ^2_5$  & 0.4205 & 0.2343 & 0.2618 &  2.1513\\                                                 $t_7$           &  0.2240 &  2.1341 &  &    $\chi ^2_7$  & 0.4358 & 0.2337 & 0.2116 &  2.1291\\
Normal          &  0.2326 &  2.0878 &  &    Lognormal    & 0.3415 & 0.2165 & 0.5487 &  2.5035\\
Logistic        &  0.2236 &  2.1429 &  &  Extreme Value  & 0.4444 & 0.2291 & 0.1872 &  2.1459
\end{tabular}
\end{center}
\end{table}

\subsection{Moments of \pdQ s}\label{sec:moments}

Let $X^*\sim F^*$, $\mu ^*=\e [X^*]$ and define the $k$th central moments by
$\mu _k^*=\e [(X^*-\mu ^*)^k]=\int _0^1(u-\mu ^*)^k\,f^*(u)\,du ,$ for $k=2,3,\dots .$
Denote the standard deviation $\sigma ^*=\sqrt {\mu _2^*}\;$ and the coefficients of skewness  and kurtosis by $\gamma _1^*= \mu _3^*/(\sigma ^*)^3$ and $\gamma _2^*=\mu _4^*/(\sigma ^*)^4$.  For all the examples in Table~\ref{table2}, $\gamma _1^*$ is a good linear predictor of the Hellinger distance of
$f^*$ from the class of symmetric distributions on [0,1], as we will see in Section~\ref{sec:closestsymm}.

\subsection{Discrete \pdQ s}\label{sec:discpdQs}
Let ${\cal X}=\{x_i\} $ be a sequence of distinct real numbers and let $\{p_i\} $ be an associated sequence of non-negative numbers $\{p_i\}$ whose sum is one. Then $\{(x_i,p_i)\}$ defines a {\em discrete probability distribution}. The class of such distributions is too rich to meaningfully discuss shape via \pdQ s, so we restrict ${\cal X}$ to a {\em lattice}, a sequence of equally spaced points. By a location-scale change this can be taken to be a subset of the integers.
For most distributions of interest to us  this domain is of the form $\{0,1,\dots ,n\}$ or $\{0,1,2, \dots \}$.  Let $\D \subset \F $ be the set of such lattice distributions; they are comprehensively treated by \cite{J-K-K-1993}.  The \pdQ  of any $F\in \D $ can now be defined in exactly the same way as for the
continuous case, but now the density function is with respect to counting measure, and the density is commonly called the {\em probability mass function} $p$ and is defined by $p(x_i)=p_i$ for all $x_i\in {\cal X}$ and $p(x)=0$ otherwise.

\begin{definition}\label{def2}
  The \cdf of each $F\in \D $ is a non-decreasing step function with jumps of size $p_i$ at $i$; it takes on value $F(x)=S_i=\sum _{j\leq i}p_j$ for $x$ lying in $[i,i+1).$ The quantile function $Q$ of $F\in \D $ is a monotone increasing left-continuous step function with value $i$ for $u$ lying in $(S_{i-1},S_i]$. The {\em discrete density quantile function} $pQ(u)$ is therefore a step function with value $p_i$ for $u$ lying in $(S_{i-1},S_i]$. In general, the graph of $pQ$ covers a sequence of adjacent {\em squares} with respective sides equal to $p_i.$  Since the sum of the areas of these squares $\kappa =\sum _ip^2_i$ is finite, we define the {\em discrete  probability density quantile} by $p^*(u)=pQ(u)/\kappa .$
\end{definition}
 Given $p^*$, one can recover $p$, but not ${\cal X}.$
 Note that a discrete distribution $p=\{ p_i\} $  on a lattice is transformed by this composition and normalization into $p^*$, which is the density with respect to Lebesgue measure of a continuous distribution with support [0,1]. It is a normalized histogram that captures the shape of the original discrete distribution on its support, free of location and scale. We often write $f^*$ for $p^*$.
  Examples of discrete \pdQ s are shown in Figures~\ref{fig2} and \ref{fig4}.

\section{Empirical \pdQ s}\label{sec:epdQ}

Given a sample of data from an unknown $F \in \F $, we want to estimate its \pdQ . This
will be done separately for the discrete  $F\in \D $ and continuous $F \in \F '$ cases.

\subsection{An empirical \pdQ \; for discrete distributions}\label{epdQdiscrete}

\begin{figure}[b!]
\begin{center}
\includegraphics[scale=.6]{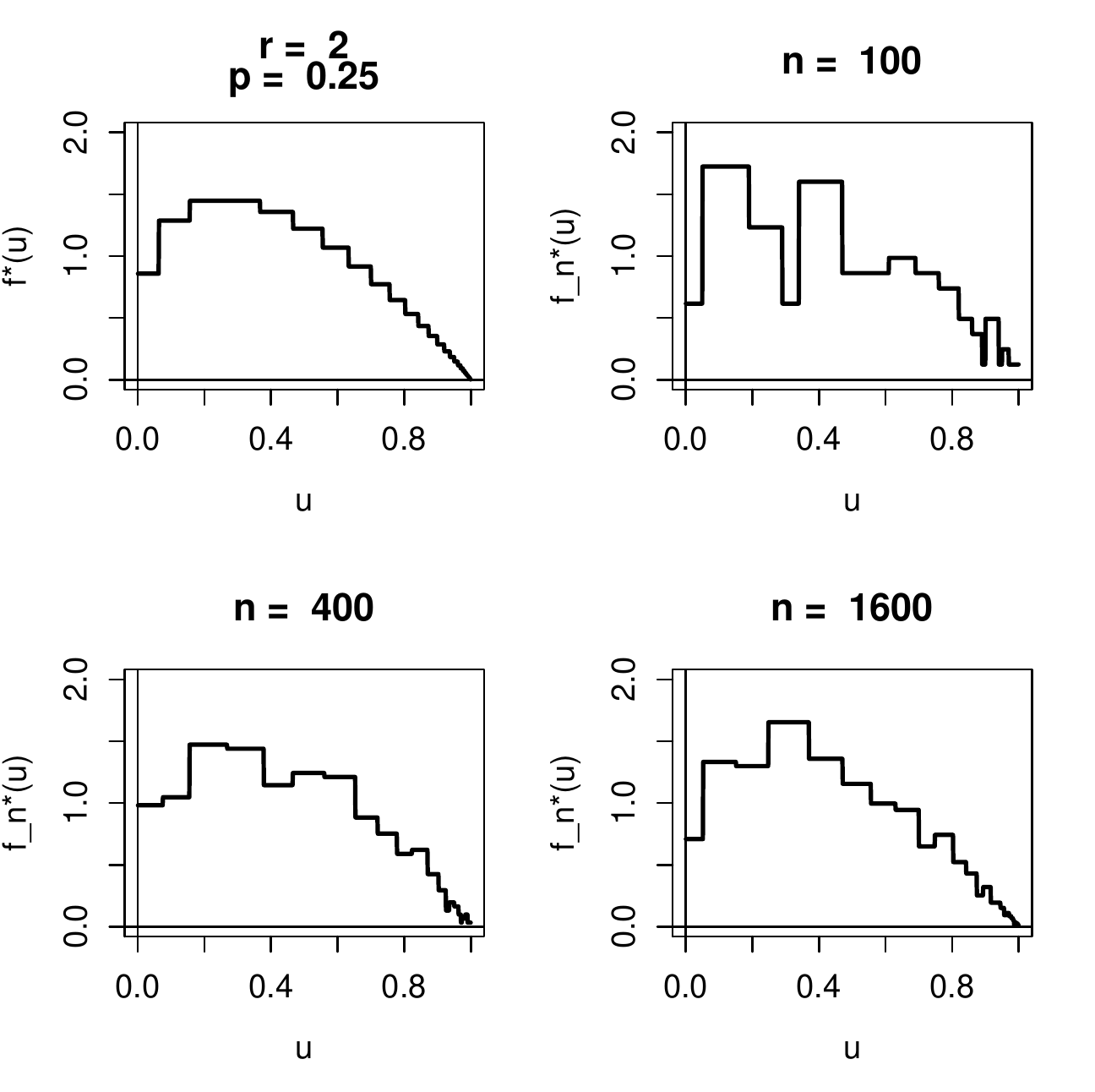}
\caption{\em The top left plot shows the graph of the \pdQ \, of  the negative binomial distribution with parameters $r=2$, $p=0.25.$ The remaining empirical \pdQ \, plots defined by (\ref{fQhatstar}) in Section~\ref{epdQdiscrete} are based on random samples of size $n$ from this distribution. \label{fig2}}
\end{center}
\end{figure}

\begin{definition}\label{def3}
Given $n$ real numbers which have $M \leq n$ distinct values $x_1 < x_2 <\dots <x_M$, let $f_n(x_m)=n_m/n$ be the relative frequency of occurrences of $x_m$ for $m=1,\dots ,M$.
Further let  $c_0=0$ and  $c_m=(\sum _{j=1}^mn_j)/n$ for $m=1,\dots ,M$.
Then the empirical \pdQ is derived in the same way as the discrete \pdQ in Definition~\ref{def2} to be:
\begin{equation}\label{fQhatstar}
    f^*_n\,(u) = \frac {n\,\sum _{m=1}^M n_m\;I\{c_{m-1} <u \leq c_{m}\}}{
    \sum  _{m=1}^M n_m^2}~.
\end{equation}
\end{definition}
If all $n$ observations are distinct, the empirical \pdQ in (\ref{fQhatstar}) is identically one, just as it is for any other discrete or continuous uniform distribution.
For another example, the top left plot of  Figure~\ref{fig2} shows the graph of the \pdQ \, of the negative binomial distribution with $r=2$ and $p=0.25.$ The other plots show empirical \pdQ s, based on random samples of varying  sample sizes. An R script for plotting the empirical \pdQ  is found in  accompanying online material.

\subsection{An empirical \pdQ \; for smooth distributions}\label{epdQcont}

For smooth distributions $F \in \F '$ we propose estimating $q(u)$ by the estimator
already studied  by many authors, including  \cite{falk-1986}, \cite{welsh-1988} and \cite{jones-1992}. It is the kernel density estimator that can be written as a linear combination of order statistics:
\begin{equation}\label{qhat}
    \hat q_n(u)=\sum _{i=1}^nX_{(i)}\,\left \{k_b\left (u-\frac {(i-1)}{n}\right )
    -k_b\left (u-\frac {i}{n}\right )\right \},
\end{equation}
where $b$ is a bandwidth and $k_b(\cdot)=\frac{1}{b}\,k(\frac {\cdot }{b})$ is the
\cite{epan-1969} kernel.

\begin{figure}[t!]
\begin{center}
\includegraphics[scale=.6]{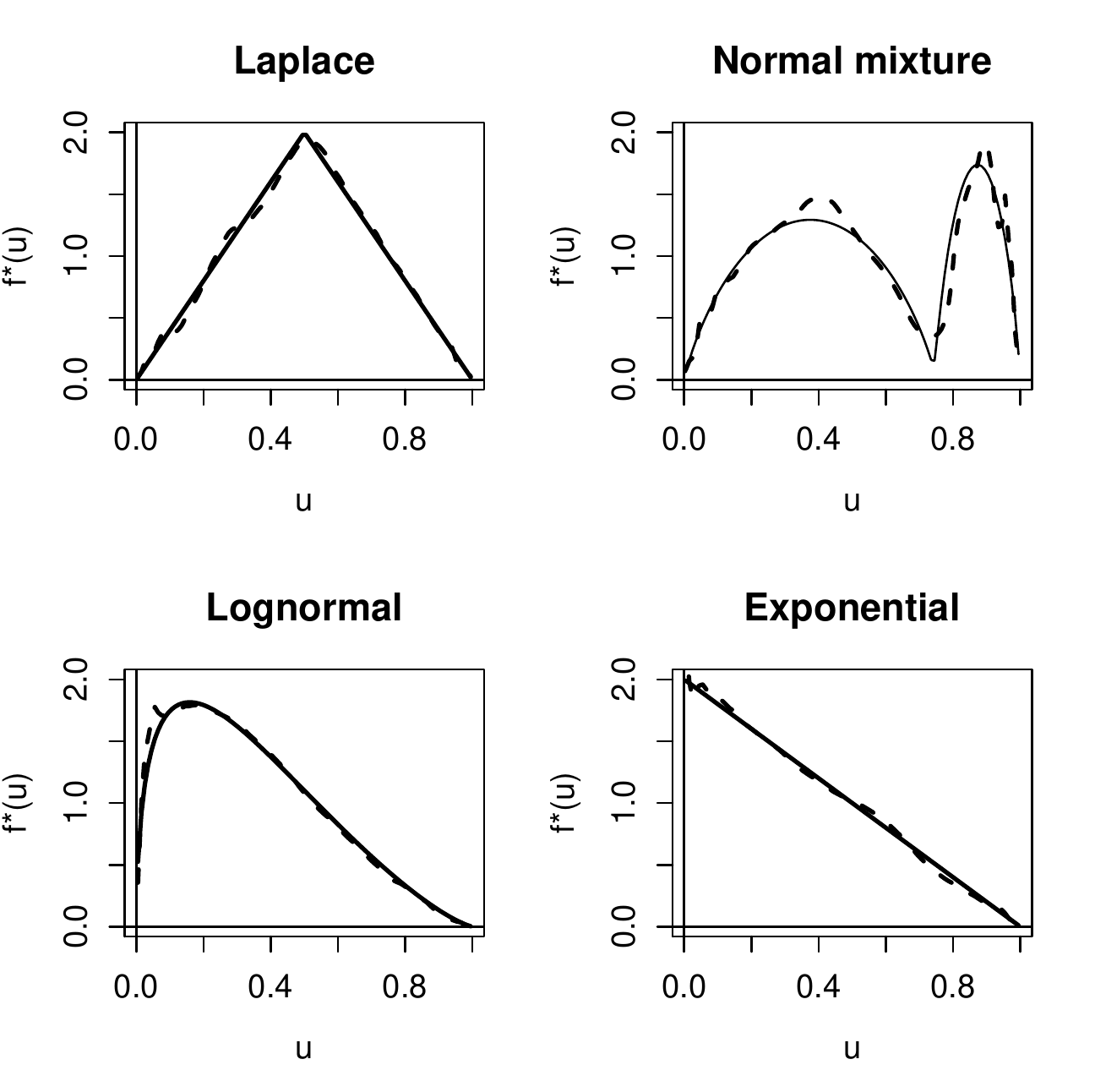}
\caption{\em  Graphs of $f^*(u)$ (solid lines) and their estimates defined  in Section~\ref{epdQcont} (dashed
lines). The normal  mixture puts weights $3/4,1/4$ on the standard normal $\Phi $ and $\Phi ((\cdot -3)/(1/4))$ respectively.  The optimal bandwidth ratio for the Cauchy is used to find the estimates in the top two plots,
while the optimal bandwidth ratio for the lognormal is the basis for the others. \label{fig3}}
\end{center}
\end{figure}

 It turns out that the  asymptotic mean squared error of $\hat q(u)$ is minimized when the bandwidth
$b(u)= (15/n)^{1/5}\; \{q(u)/q''(u)\}^{2/5}.$ The ratio $q(u)/q''(u)$ is similar in shape to the density quantile $fQ(u)=1/q(u)$, and hence remarkably stable for $F$ in broad classes such as all symmetric unimodal distributions, or all $F$ with
positive unimodal density on $[0,+\infty ).$ \cite{prst-2016} show that by employing the optimal ratio for the Cauchy, one obtains good estimators $\hat q(u)$ of $q(u)$ for all $F$ in the first
class, while the optimal ratio for the lognormal yields good estimators for  all $F$ in the
second.
In Figure~\ref{fig3} are shown examples of the graphs of $f^*(u)=1/\{\kappa q(u)\}$ and the
estimators $f_n^*(u)=1/\{\hat\kappa \hat q_n(u)\}$, $u=0.005:0.995/0.01$, for sample size $n=1600$.
Sample sizes as low as 400 give a good idea of the shape.

\section{Applications}\label{sec:applics}

We find location- and scale-free global distances between shape families of distributions, which
leads to effective fitting of shape families to data.

\subsection{Divergence and distance measures}\label{metrics}

Given densities $f_1$, $f_2$ with respect to Lebesgue measure, the {\em Hellinger distance}  between them is defined by $ H(f_1,f_2)= [2^{-1}\,\int \{\sqrt {f_1(x)}\, -\sqrt {f_2(x)}\,\}^2\;dx ]^{1/2}.$
The {\em Kullback-Leibler information} $I(f_1:f_2)$ in $X\sim f_1$ for discrimination between $f_1,f_2$
is defined by \cite[p.5]{Kullback-1968} as $I(f_1:f_2) = \int \ln (f_1(x)/f_2(x))\; f_1(x)\, dx.$
The {\em symmetrized divergence}, or \kld, is defined by
$J(1,2)=I(1:2)+I(2:1).$ We often abbreviate $H(f_1,f_2)$ to $H(1,2)$, $I(f_1:f_2) $ to $I(1:2)$ and $J(f_1,f_2)$ to $J(1,2)$. Further, we denote by $H^*(1,2)$ the Hellinger metric applied to the \pdQ s $f_1^*, f_2^*$ of $f_1,f_2,$ and similarly for $I^*(1:2)$ and $J^*(1,2).$
 An advantage of working with $f^*$s over $f$s is that the \pdQ s for both discrete and continuous $F$ are densities with respect  to Lebesgue measure, so the Hellinger distance or \kld between their \pdQ s can be informative.

\subsection{Probabilistic examples}

\subsubsection*{Ex.~1:\  Rate of convergence in Hellinger distance of Poisson from normality.}
It is well known \cite[p.161]{J-K-K-1993} that  $X\sim $Poisson($\lambda $)  approaches the normal distribution as $\lambda \to \infty $ in that $(X-\lambda )/\sqrt {\lambda }\to Z\sim \Phi $ in distribution. And numerical computations show that the  Hellinger distance between the Normal \pdQ and that of the Poisson($\lambda )$ distribution is $H=H(f_\lambda ^*,\varphi ^*) \approx 0.17077/\sqrt {2.4\, \lambda -1}\,$ for $1\leq \lambda \leq  100,000.$

\subsubsection*{Ex.~2:\  Discretized exponential.}

Given $X\sim \text{Exponential}(1$), for $0<r\leq 1$ let $Y_r=\lfloor rX\rfloor \sim $ Geometric($p_r$), with $p_r=1-e^{-r}.$ Then the Hellinger distance of $X^*$ from $Y_r^*$ is $H_r\approx r/10$; and the root-KLD divergence $\sqrt {J_r}\approx 3r/11,$ for $0<r\leq 1.$ This example illustrates that one can quantify the Hellinger distance between a continuous distribution and a discrete approximation to it.

\subsubsection*{Ex.~3:\  Matching symmetric Tukey distributions with Student's t-distributions.}

Student's t-distributions are well approximated by certain symmetric Tukey($\lambda $) distributions, see \cite{tukey-1960}, \cite{join-1971}, \cite{rogers-1972} and references therein. These authors used expected ranges or specific tail probabilities to match distributions; and in particular find that for the
Cauchy $(\nu=1)$ the symmetric Tukey distribution with $\lambda =-1$ is close to it, while for the Normal $(\nu \to \infty )$ the Tukey with $\lambda = 0.14$ is quite close. Here we find the best global fits using metrics on the respective \pdQ s, and obtain very similar results. In Table~\ref{table3} are listed the $\lambda =\lambda_\text{min}(\nu ) $ that minimize the Hellinger distance between the \pdQ s of the Tukey($\lambda $) family and the \pdQ of the Student's t distribution with $\nu $ \df . The minimum distance $H_{\min}$ is also shown. For the normal distribution we found $\lambda _{\min}(\infty )=0.14435$ with $H_{\min}=0.0010.$  A good approximation for $\nu \geq 12$ is given by $\lambda _{\min}=0.14435-1/(1.07\,\nu).$  Minimizing the KLD instead of
the Hellinger metric led to the same results for $\lambda _{\min}.$
\begin{table}[t!]
\caption{\small\em \label{table3} For selected values of $\nu$ are shown the  $\lambda _\text{min}$
 to 3 decimal places that minimizes the Hellinger distance of the \pdQ s for the symmetric Tukey($\lambda $) distributions from the \pdQ \ of Student's t-distribution with
 $\nu $ degrees of freedom. Also shown are the minimum Hellinger distances.}
 \begin{footnotesize}
 \begin{center}
 \begin{tabular}{crrrrrrrrrrr}
 $\nu $ & 0.5 & 1 &  2 & 3& 4 & 5& 6 & 7 &  12 & 24 & 100   \\
 \hline
 $\lambda _{\min}$ & $ -1.881$ & $-0.867$ & $-0.357$ & $-0.188$ & $-0.104$ & $-0.053$ & $-0.020$ & 0.004 &  0.063 & 0.104 & 0.135 \\
 $H_{\min}$ & 0.005 & 0.005 & 0.004 & 0.003 & 0.003 & 0.003 & 0.002 & 0.002 & 0.002 & 0.001 & 0 .001
 \end{tabular}
 \end{center}
 \end{footnotesize}
 \end{table}

\subsection{Data examples}
 Robust global fitting of shape parameter families to data is possible by minimizing the Hellinger distance
of the model \pdQ\  to the empirical \pdQ , as we now demonstrate for three important shape families.

\subsubsection*{Ex.~4:\  Fitting symmetric Tukey models to data.}

A standard method of fitting a family of distributions indexed by a shape parameter to data is due to \cite{fill-1975}.
He suggested finding the correlation coefficient between the ordered data and the quantiles determined by the family for a range of  values of the shape
parameter and then constructing the \lq probability plot correlation coefficient\rq\  (ppcc): \ the correlation coefficient as a function of shape parameter.  The shape parameter that maximizes the coefficient  is his proposed ppcc estimate.
This method is often applied to symmetric-looking data to fit a member of the Tukey($\lambda $) family, and if an estimate of $\lambda $ were near $-1$, (or 0.14), say, it is suggested to assume a Cauchy (or normal) model, see the discussion in  Example~3 and Table~\ref{table3}. Location and scale estimates are then obtained by regressing the sorted data on the quantiles of the chosen shape model.

As an alternative method of fitting a family of distributions indexed by a shape parameter to data, we propose finding the empirical \pdQ \ $f_n^*(u)=1/\{\hat\kappa \hat q_n(u)\}$ of Section~\ref{epdQcont} for $u$ in a grid on
[0,1]. This nonparametric density estimate requires a sample size of at least 500.
Then, for each shape parameter $\lambda $, say, in a large range of values, compute the Hellinger distance $H_\lambda =H(f_n^*,f_\lambda ^*)$ of $f_n^*$ from the model \pdQ \;$f_\lambda ^*$. The $\lambda $ that minimizes this distance will be called the $H$-\pdQ \; {\em goodness-of-fit estimate} of $\lambda .$

Both methods require a search over a grid of $\lambda $ values and we recommend that this be done in two stages:
first, for a rough grid over a large range, say increments of size 0.2 over $[-10, 10]$; and second, over a narrow range around the first estimate with increments of size 0.01, say.   In our simulation study we made one pass with increments 0.01 over $[-10, 10]$.  For sample size 500 the run time for this search was approximately 0.6 second.
To compare these methods we ran 400 replications of an experiment
with sample sizes 500 selected at random from the symmetric Tukey($\lambda $) distribution using the R package {\tt gld} due to \cite{gld} with the \cite{RS} parameterization. The R command to obtain a sample of simulated
data $x$ is {\tt x <- rgl(500, c(0,lambda,lambda,lambda), param = "rs", lambda5 = NULL)}.
The choices of $\lambda $ listed in Table~\ref{table4} are representative of very long, long, normal and truncated tails, respectively, also discussed later in Section~\ref{sec:tailwt}.

The results of fitting these data sets by the methods ppcc and $H$-\pdQ \; are summarized in Table~\ref{table4}. In particular, we list the empirical standard errors (SE equals the square root of the sample  variance plus squared estimated bias).
The rows corresponding to  $\lambda =-1$ and $\lambda =-2$ suggests that for long and very long tails, the \pdQ approach is much more efficient at identifying the data source.  However, for the cases where the data are approximately normal or have truncated tails, there is little to choose between the methods.

\begin{table}[h!]
\caption{\small\em \label{table4} For selected values of $\lambda $ are shown summary results for 25
competing estimates (ppcc and $H$-\pdQ) based on samples of size 500 from the Tukey($\lambda $).}
 \begin{footnotesize}
 \begin{center}
 \begin{tabular}{crrrrrrr}
 $\lambda $& Method & SE   &   & mean       & sd      & min         & max      \\[.2cm]
\hline
 $-2$  & ppcc       & 2.61 &   &  $ -3.31$  & 2.26    & $-10.00 $    & $ -0.86$ \\
       & $H$-\pdQ   & 0.18 &   &  $ -2.02$  & 0.18    & $-2.76 $    & $-1.52 $ \\
 \cline{2-8}
 $-1$  & ppcc       & 1.46 &   &  $ -1.72$  & 1.27    & $-8.41$     & $ -0.63$ \\
       &  $H$-\pdQ  & 0.11 &   &  $ -1.02$  & 0.11    & $-1.34 $    &  $-0.68$ \\
 \cline{2-8}
 $0.14$& ppcc       & 0.05 &   &    0.12    & 0.05    & 0.01        &  0.28     \\
       &  $H$-\pdQ  & 0.06 &   &    0.14    & 0.06    & 0.05        &  0.29     \\
 \cline{2-8}
 $3  $ &  ppcc      & 0.94 &   &    2.67    & 0.88    & 0.45        &   3.52    \\
       &  $H$-\pdQ  & 0.60 &   &    2.77    & 0.55    & 0.52        &   3.43
\end{tabular}
 \end{center}
 \end{footnotesize} \end{table}

\subsubsection*{Ex.~5:\ Fitting Weibull models to data.}

Another shape family that is often fitted to lifetime or income data is the Weibull model with shape  parameter $\beta $. We compare three methods: the ppcc and $H$-\pdQ \
 methods just described in the last example; the third is the maximum likelihood (MLE) approach, available through the command {\tt fitdistr(x,"weibull")} on the software \cite[Core Development Team]{R}  with the internal package MASS, \cite{MASS}; it assumes $x>0$ and returns both scale and shape MLEs and their standard errors based on the
 observed information matrix.

The summary results for 400 replications of sample size 500 from six data configurations are listed in
Table~\ref{table5}. For the Weibull  $\beta =1$  (exponential) model the MLE approach performs best, having the smallest SE, while the $H$-\pdQ \  method is the second best performer.
In the second and third configurations of 5\% contaminated data,
  the MLE and $H$-\pdQ \  methods are comparable, while the ppcc method is badly affected by a few outliers.
The simulations for $\beta =2$ show that the MLE approach can also perform poorly for contaminated data.

Similar results to those above  were obtained by fitting the Gamma shape family.

\begin{table}[h!]
\caption{\small\em \label{table5}     Summary results of
competing estimates (ppcc, $H$-\pdQ and MLE) of the Weibull shape parameter $\beta $  based on 400 replications of
samples of size 500 from the Weibull $W_\beta $ with $\beta =1,2$ distributions; and also from 95:05
mixture distributions of $W_\beta $ with contamination from the standard lognormal {\em (LN)} and then such contamination multiplied by two {\em(2\,LN)},}
 \begin{footnotesize}
 \begin{center}
 \begin{tabular}{lrrrrrrr}
 Data Source      & Method         & SE    &&  mean & sd   & min   & max   \\[.2cm]
\hline
                  & ppcc           & 0.09  &&  0.96 & 0.09 &  0.66 & 1.25    \\
 $W_1  $          & $H$-\pdQ       & 0.06  &&  0.97 & 0.06 &  0.64 & 1.09    \\
                  & MLE            & 0.04  &&  1.00 & 0.04 &  0.89 & 1.12    \\
  \cline{2-8}
                  & ppcc           & 0.17  &&  0.89 & 0.14 &  0.22 & 1.17    \\
 $0.95\,W_1 +0.05(LN)$& $H$-\pdQ   & 0.05  &&  0.98 & 0.05 &  0.81 & 1.14    \\
                  & MLE            & 0.04  &&  1.00 & 0.04 &  0.88 & 1.10    \\
  \cline{2-8}
                    & ppcc         & 0.42  &&  0.64 & 0.21 &  0.10 & 1.09    \\
$0.95\,W_1 +0.05(2\,LN)$& $H$-\pdQ & 0.06  &&  0.97 & 0.05 &  0.79 & 1.10    \\
                    & MLE          & 0.07  &&  0.94 & 0.05 &  0.79 & 1.04    \\
  \hline
                  & ppcc           & 0.15  &&  1.99 & 0.15 &  1.56 & 2.47    \\
  $W_2  $         & $H$-\pdQ       & 0.18  &&  2.12 & 0.14 &  1.71 & 2.52    \\
                  & MLE            & 0.07  &&  1.96 & 0.07 &  1.78 & 2.23    \\
  \cline{2-8}
                  & ppcc           & 1.12  &&  0.99 & 0.48 &  0.10 & 2.13    \\
$0.95\,W_2+0.05(LN)$ & $H$-\pdQ    & 0.12  &&  1.97 & 0.12 &  1.66 & 2.40    \\
                  &  MLE           & 0.39  &&  1.65 & 0.19 &  1.01 & 2.10    \\
  \cline{2-8}
                  & ppcc           & 1.61  &&  0.40 & 0.22 &  0.10 & 1.29    \\
$0.95\,W_2 +0.05(2\,LN)$& $H$-\pdQ & 0.15  &&  1.89 & 0.11 &  1.50 & 2.16    \\
                  &  MLE           & 0.68  &&  1.34 & 0.16 &  0.95 & 1.80
  \end{tabular}
 \end{center}
 \end{footnotesize}
 \end{table}

\subsubsection*{Ex.~6:\ Fitting shape models to wool fibre diameter data.}

Raw wool is routinely classified by the distribution of fibre diameters, measured in microns. A
typical example of 4817 measurements obtained by laser scanning technology is given in Appendix~\ref{sec:appenda}. Of main interest to wool assessors is the mean diameter $\bar x =25.08$, the standard deviation $s=5.388$ and the coefficient of variation $cv=0.215,$ as well as percentages of
relatively large fibres that can cause prickliness in finished woolen goods. Here we fit
parametric models to these data using the three methods described above.

\begin{table}[h!]
\caption{\small\em \label{table6} The best fitting  Weibull($\beta $)
and Gamma($\alpha $) models to the wool fibre diameter data are listed,
 along with location and scale estimates obtained by regressing the quantiles of
 the sorted data on the quantiles of the respective fitted shape models. The values of $H$ are Hellinger distances of the empirical \pdQ \ of Section~\ref{sec:epdQ} from each of the \pdQ s with estimated shape parameters.}
 \begin{footnotesize}
 \begin{center}
 \begin{tabular}{lcrrrrc}
 Family  & Method     &  Shape           & Location    & Scale         && $H$      \\[.2cm]
\hline
         & ppcc       &  2.44            &  12.80      & 13.85         && 0.0360       \\
 Weibull & $H$-\pdQ   &  2.76            &  11.38      & 15.40         && 0.0323       \\
         & MLE        &  4.82            &   2.66      & 24.47         && 0.0595       \\
   \cline{2-7}
         & ppcc       &  22.21           &  $-0.26$    & 1.14          && 0.0234       \\
 Gamma   & $H$-\pdQ   &  35.75           &  $-7.09$    & 0.89          && 0.0220       \\
         & MLE        &  21.67           &  0.03       & 1.15          && 0.0236       \\
   \end{tabular}
 \end{center}
 \end{footnotesize}
 \end{table}
For the Weibull family, the  ppcc and $H$-\pdQ \, results are very similar.
 The MLE estimate of $\beta $ is quite different, because it was obtained by the MLE optimization method that assumes zero location and estimates scale and shape simultaneously. That is why the regression results for location and
scale shown in Table~\ref{table6}, namely 2.66 and 24.47 are close to zero and the MLE estimate of scale 27.26,
respectively. The Hellinger distance of the MLE fitted \pdQ\, from the empirical \pdQ \
is almost twice that of the distances for the other two fitted models.

For the Gamma models, the Hellinger distances for all three methods are all about 1/2 those for the Weibull models. QQ-plots for these six best fitting distributions reveal that the Weibull models all fit quite badly; while for the
fitted gamma models all QQ-plots are almost linear. Examples of some QQ-plots are in Figure~\ref{fig6} of Appendix~\ref{sec:appenda}.  For the gamma family, a wide range of shape parameters leads to a good fit; and an
examination of the ppcc plot reveals it is almost flat for $20 <\alpha <50.$  Also, the corresponding Hellinger distances are almost the same for this range of $\alpha $.

\section{Measuring asymmetry in a \pdQ }\label{sec:asymmetry}

In addition to skewness and kurtosis of \pdQ s, the \lq shape\rq\  can be described by the degree of asymmetry, the location of modes, and the tail behavior of $f^*$ at 0 and 1.

\subsection{Finding the closest symmetric distribution to a \pdQ }\label{sec:closestsymm}

\begin{figure}[b!]
\begin{center}
\includegraphics[scale=.6]{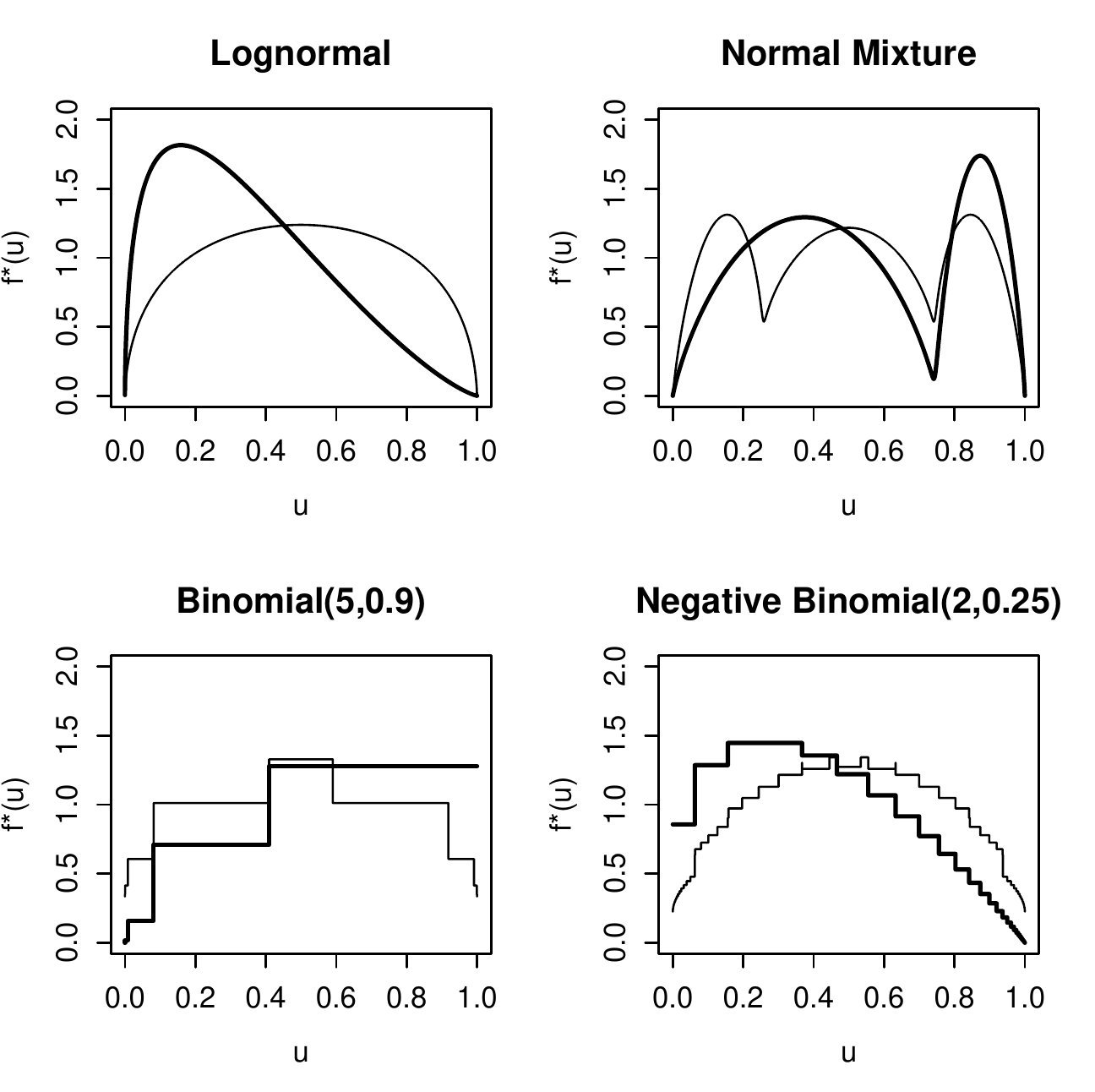}
\caption{\footnotesize  \em  Graphs of $f^*(u)$ in thick black lines and their respective closest (in the
Hellinger distance) symmetric distributions $f^*_\text{symm}$ in thin lines. The normal  mixture puts weights $3/4,1/4$ on the standard normal $\Phi $ and $\Phi ((\cdot -3)/(1/4))$ respectively. The respective Hellinger distances to the symmmetric class on [0,1] in the four cases are  0.2386, 0.1255, 0.1376 and  0.1443.\label{fig4}}
\end{center}
\end{figure}

One can measure the skewness of a \pdQ with the coefficient $\gamma _1^*$ as in Table~\ref{table1},
or with many other definitions found in \cite{S-2014}. However, asymmetry entails more than skewness,
because asymmetric distributions can have 0 skewness.
Here we measure the asymmetry of a \pdQ by finding its distance or divergence to the class of symmetric distributions on [0,1].
It is shown in Theorem 3.2 of \cite{with-2010} that if $f_1,f_2$ are densities with respect to Lebesgue measure
and $f_2$ is assumed to be symmetric about 0, then the closest such $f_2$ to $f_1$ in terms of minimizing the Hellinger distance is given by $f_2(x)=\alpha ^2(x)/d$ where $\alpha (x)=\{\sqrt {f_1(x)}\,+\sqrt {f_1(-x)}\,\}/2$ and
$d=2\int _0^{+\infty}\alpha ^2(x)\;dx.$
Consequently, if $f^*$ is an arbitrary \pdQ , the Hellinger-closest symmetric distribution is given by
 $f^*_\text{symm}(u)=\{\alpha ^*(u)\}^2/d^*$ where $\alpha ^*(u)=\{\sqrt {f^*(u)}\,+\sqrt {f^*(1-u)}\,\}/2$ and $d^*=\int _0^1\{\alpha ^*(u)\}^2\;du.$  The actual Hellinger distance $H_\text{min}(f^*)=H(f^*,f^*_\text{symm})$ of $f^*$ from the class of symmetric distributions is determined by:
\begin{equation}\label{Hdisttosymm}
     2\{1- H_\text{min}^2(f^*)\}^2= 1+\int _0 ^1\sqrt {f^*(u)f^*(1-u)}\;du  ~.
\end{equation}
Some examples are shown in Figure~\ref{fig4}.

\cite{with-2010} also find the closest symmetric distribution to a given one in the sense of minimizing the Kullback-Leibler
divergence $I(2:1).$ Their result is extended to minimizing  $J(1,2)$ in the following proposition.

\begin{proposition}\label{prop1}
Given $f_1,f_2$ probability densities with respect to Lebesgue measure, and assume that $f_2$ is symmetric about 0. Then
\begin{description}
  \item[(a)] The $f_2$ that minimizes $I(2:1)$ is given by $f_2(x)=\nu (x)/d$, where $\nu (x)= \sqrt {f_1(x)\,f_1(-x)}\,$ and $d=2\int _0^{+\infty}\nu (x)\;dx.$
  \item[(b)] The $f_2$ that minimizes $I(1:2)$ is given by $f_2(x)\equiv \overline f_1(x)\equiv \{f_1(x)+f_1(-x)\}/2$.
  \item[(c)] The $f_2$ that minimizes $J(1,2)$ is given by $f_2(x)=\beta (x)/d$ where $\beta (x)$ is the solution to:
\begin{equation}\label{Jmin}
     \beta (x) = \nu (x) \,\exp \left \{\overline f_1(x)/\beta (x)\right \}~;\\
      \end{equation}
 and $d=2\int _0^{\infty }\beta (x)\; dx .$
\end{description}
\end{proposition}

The proof of Proposition~\ref{prop1} is in the Appendix~\ref{sec:appendb}, along with an algorithm for computing the minimizer of $J$ guaranteed by Part (c) and two illustrative examples.

In our applications of Proposition~\ref{prop1} to \pdQ s $f_1^*, f_2^*$ on
the unit interval, where symmetry takes place about 1/2,  $\overline f_1(x)\equiv \{f_1(x)+f_1(-x)\}/2$ becomes $ \overline f_1^*(u)\equiv \{f_1^*(u)+f_1^*(1-u)\}/2$, and similarly  $\nu(x)$ becomes $\nu ^*(u)\equiv \sqrt {f_1^*(u)\,f_1^*(1-u)}\,$.  The minimum divergences realized in applying
parts (a)-(c) of Proposition~\ref{prop1} to the \pdQ s are denoted, respectively, $I^*_{\ :1}$, $I^*_{1:\ }$
and $J^*$.

\begin{table}[h!]
\begin{small}
\begin{center}
\caption{\label{table7}\em \small
The coefficient of skewness $\gamma _1 ^*$ of $f^*$ taken from
Table~\ref{table2} is listed, along with the  $H^*=H_\text{min}(f^*)$ of (\ref{Hdisttosymm}), and  the minimum divergences $I^*_{1:\ }$, $I^*_{\ :1}$ and $J^*$ of Proposition~\ref{prop1} applied to $f^*_1=f^*$.}
\begin{tabular}{cccccc}
\hline
$F$     & $\gamma _1 ^*$ & $H^*$  & $I^*_{1:\ }$ & $I^*_{\ :1}$ & $J^*$ \\
\hline
    Pareto(0.5)  & 1.0498 & 0.4421 & 0.5401 & 1.2224 & 2.1589  \\
    Pareto(1)    & 0.8607 & 0.3660 & 0.4077 & 0.6931 & 1.2710  \\
    Pareto(2)    & 0.7318 & 0.3094 & 0.3107 & 0.4535 & 0.8507  \\
   Weibull(2)    & 0.1315 & 0.0672 & 0.0178 & 0.0182 & 0.0363  \\
   $\chi ^2_2$   & 0.5657 & 0.2349 & 0.1931 & 0.2416 & 0.4646  \\
   $\chi ^2_3$   & 0.3808 & 0.1687 & 0.1061 & 0.1191 & 0.2326  \\
   $\chi ^2_5$   & 0.2618 & 0.1191 & 0.0548 & 0.0580 & 0.1145  \\
   $\chi ^2_7$   & 0.2116 & 0.0970 & 0.0368 & 0.0382 & 0.0757  \\
   Lognormal     & 0.5487 & 0.2386 & 0.2014 & 0.2500 & 0.4747  \\
 Extreme Value   & 0.1872 & 0.0855 & 0.0287 & 0.0295 & 0.0587
\end{tabular}
\end{center}
\end{small}
\end{table}

\begin{figure}
\begin{center}
\vspace{-1.5cm}
\includegraphics[scale=.6]{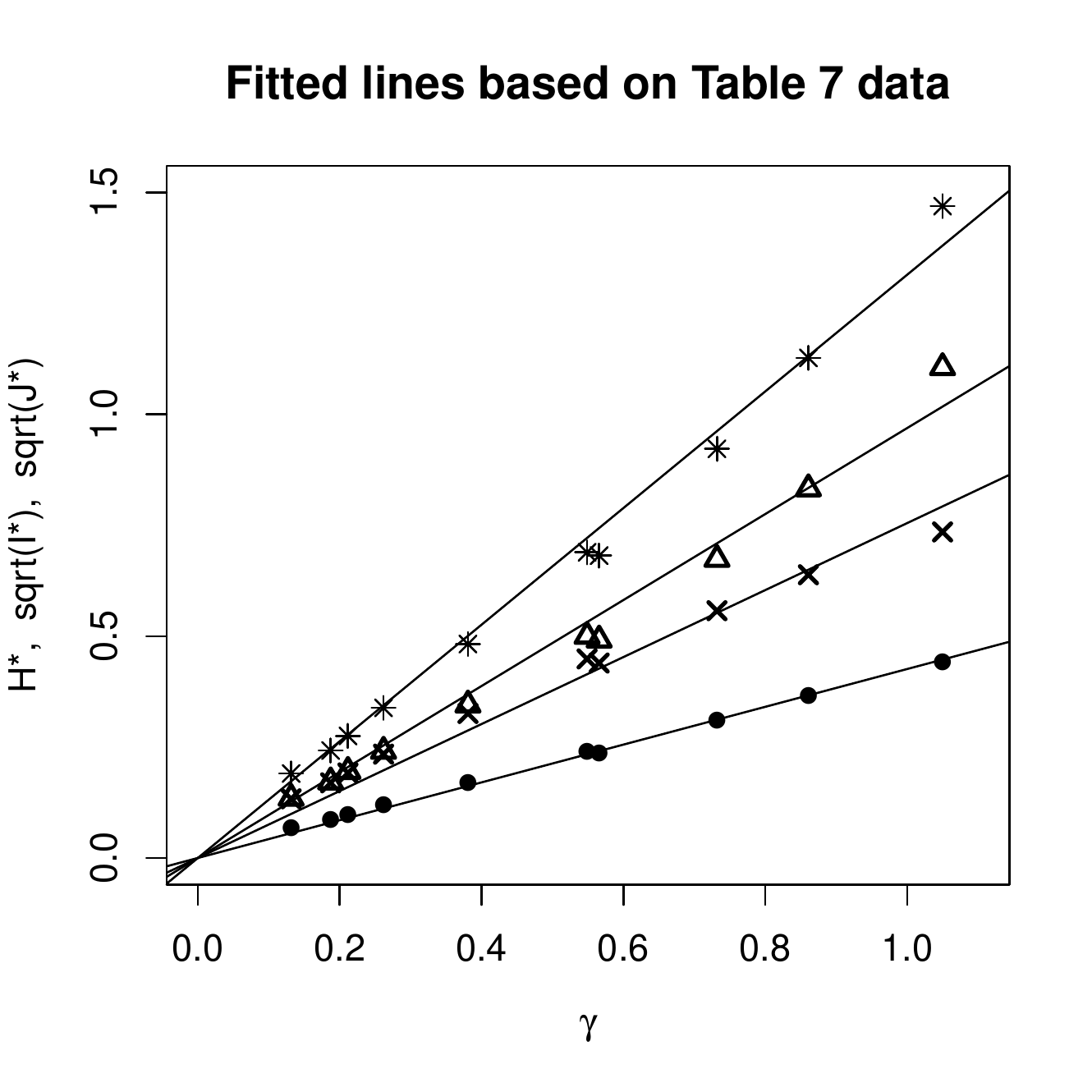}
\caption{\footnotesize  \em  Plots of $H^*$ versus  $\gamma \equiv \gamma _1 ^*$ in circles, and $\sqrt {I^*_{1:}}\,$, $\sqrt {I^*_{:1}}\,$ and  $\sqrt {J^*}\,$ against  $\gamma $ in crosses, triangles and asterisks, respectively,
based on results in Table~\ref{table7}. Superimposed are their least-squares lines through the origin, with
approximate coefficients 0.43, 0.75, 0.97 and 1.31. \label{fig5}}
\end{center}
\end{figure}

In Table~\ref{table7} are listed the values of some possible measures of asymmetry for a given \pdQ \ $f^*$ for $f$
belonging to an asymmetric location-scale family $F$.
These values  are positively  correlated, and
$\gamma _1 ^*$, $H^*$,  $\sqrt {I^*_{1:\ }}\,$, $\sqrt {I^*_{\ :1}}\,$ and  $\sqrt {J^*}\,$ are very highly correlated, as shown in Figure~\ref{fig5}.

To show that $H^*=0.43\cdot |\gamma _2^*|$ is only an approximation on the class of \pdQ s, consider the power function family $F_b(u)=u^b,$ $0<u<1$,  $b>0$, which is also the Beta($b,1$) family.  Assuming $b>1/2$, it has \pdQ \ $f^*_b=f_{b^*}$ with $b^*=2-1/b,$ so $f_b^*\sim $Power$(b^*),$ or Beta($b^*,1$). The skewness coefficient
of the Beta family is known \cite[p.217]{J-K-B-1995}, and in this case  $\gamma _1(b^*)=2(1-b^*)\sqrt{1+2/b^*}\,/(b^*+3),$
for $\beta ^*\in (0,2]$. It is monotone decreasing in $b^*$ over this range with  limiting values
$\gamma _1(0+)=+\infty $  and $\gamma _1(2)=-2\sqrt 2\;/5.$

 To find the Hellinger distance $H_\text{min}(f_{b^*})$ from the symmetric class, we need \[\int _0^1\sqrt {f_b^*(u)f_b^*(1-u)}\,\;du =b^* \text{Beta}\left (\frac {b^*+1}{2},\frac {b^*+1}{2}\right )=\frac{\Gamma ^2\left(\frac {b^*+1}{2}\right)}{\Gamma (b^*)}~.\]
This quantity approaches $\pi /4$ as $b^*\to 2$, and using (\ref{Hdisttosymm}) the minimum Hellinger distance to a symmetric distribution approaches   $\sqrt {1-\sqrt {(1+\pi/4)/2 }\;}= 0.2348882$. The ratio $|\gamma _1(b^*)|/H_\text{min}(f_{b^*})\to  0.415228$ as $b^*\to 2,$ a value not far from the  gradient 0.43 fitted to the 10 models in Figure~\ref{fig5}. The ratio
$|\gamma _1(b^*)|/H_\text{min}(f_{b^*})$ is indeed greater than 0.37 for $b^*>1$,  but for $b^*\leq 2/3$ this ratio falls below 1/3. This example suggests that the approximation $H_\text{min}(f^*)\approx 0.43\,| \gamma _1(f^*)|$ may only hold for bounded $f^*.$

\section{Classification by tail-weight}\label{sec:tailwt}

We focus on the right-hand tails, leaving the  adjustments for left-hand tails to the reader.
Recall that $f^*(u)=fQ(u)/\kappa $ so $f^*(u)$ is the rate at which the density $f$ is accumulating
probability at its $u$th quantile. Given $f^*,g^*$ with $f^*(u)/g^*(u)>1$
for $u$ approaching 1, we would say $f$ has the {\em shorter} right tail, because it reaches its large quantiles faster than $g$ reaches its corresponding large quantiles. This enables one to partially order tails of all densities $f$ with support bounded on the right, which is perhaps the least interesting case.
We somewhat arbitrarily define  such $f$ to have \lq short\rq\  right tails.

The more interesting cases are for $f$ with unbounded support to the right, in which case $f^*(1)=0.$
The first  derivative $(f^*)'(u)=f'(Q(u))\,q(u)=-J(u)/\kappa $, where
$J(u)=- (f'/f)(Q(u))$ is the score function for location-scale families that plays an important role
in classical nonparametric statistics, see e.g. \cite{hajek-1967}.
Now, because  $f^*(1)=0$,
\begin{equation}\label{fstarpr1}
 (f^*)'(1)\; =\; \lim _{u\to 1} \frac{f^*(1)-f^*(u)}{1-u}\;=\;
\lim _{u\to 1}\,\frac{-fQ(u)}{\kappa (1-u)}\;\leq \;0~.
\end{equation}
\cite{par-1979} observes that when $fQ(u)/(1-u)^\alpha $ approaches a finite positive limit for some
positive $\alpha $, then the intervals  $0<\alpha <1$, $\alpha =1$ and $\alpha >1$ correspond
to the statistician's perception that  probability laws have three types of behavior: short, medium (exponential) and long, respectively. Note that $(f^*)'(1)$ in (\ref{fstarpr1}) is essentially case $\alpha =1.$

\subsection{Definitions and properties}\label{sec:tailwtdefns}

The above considerations plus examination of several examples below
lead us to introduce a more detailed description of tail behavior, based on
derivatives of the \pdQ at the boundaries. This  system not only  divides tails
into classes based on their relative rate of convergence to 0 relative to the exponential, but also seeks within classes to provide an absolute measure of tail-weight within classes.

\begin{definition}\label{def4}
For $n\geq 1$ define
$(f^*)^{(n)}(1)\equiv \lim _{u\uparrow 1}(f^*)^{(n)}(u)$, assuming  these derivatives exist as finite numbers for $u$ near one and their limits as $u\to 1$ exist as finite or infinite values.

 Let $n^*$ be the smallest integer $n\geq 1$ for which
which $(f^*)^{(n)}(1)\neq 0$. If $n^*$ exists, the right tail is of $n^*$-order; otherwise it
is of infinite $*$-order.  The  tail is called  medium, long or very long respectively, if $n^*=1,2$ or $n^*\geq 3$.
\end{definition}

 A sampling of examples is  in  Table~\ref{table8}.
Note that for medium tails (the case $n^*=1$), one has   $(f^*)(1)=0$ and  $(f^*)'(1)\neq 0$, but this derivative at 1 must be non-positive.  The larger the magnitude $|(f^*)'(1)|,$ the shorter the tail. Similarly for long-tailed \pdQ s, the case $n^*=2$, one has   $(f^*)(1)=(f^*)'(1)=0$ and  $(f^*)''(1)\neq 0$. Now $(f^*)'(u)\leq 0$ for $u$ near 1, so its derivative $(f^*)''(1)$ will be non-negative or $+\infty $.

The \pdQ s are not always ordered because different \pdQ s, such as the normal and Tukey with $0<\lambda <1$, both have infinity limits of their derivatives as $u\to 1$. One could order such \pdQ s by comparing the ratio of their derivatives $(f^*)'(u)/(g^*)'(u)$ as $u\to 1$, and similarly for distributions with long or very long tails, but we will not pursue this here.

\subsection{Examples}\label{tailexs}

\subsubsection*{The lognormal distribution.}
 The lognormal \pdQ is from Table~\ref{table1} given by
 $f^*(u)=c\,\varphi(z_u)\exp (-z_u),$ where $z_u=\Phi ^{-1}(u)$ and
 $c=1/\kappa =2\sqrt {\pi }\;\exp (-0.25).$  Thus $(f^*)'(u)=-c\,\exp (-z_u)(1+z_u)$, which approaches 0 as $u \to 1.$   Further, $(f^*)''(u)=c\,z_u\,\exp (-z_u)/\varphi(z_u)$, which approaches $+\infty $ as $u\to 1$ so the lognormal
 right tail is \lq long\rq .\

\begin{table}[t!]
\begin{small}
\begin{center}
\caption{\label{table8}{\bf Right-tail behavior
 of some continuous distributions.} \em  The \pdQ s are listed in Table~\ref{table1}. Note that if $f^*(1)=0$, then  $(f^*)'(1)\leq 0$.}
\vspace{.25cm}
\begin{tabular}{llcrcc}
 tail  &     $F$                               && $(f^*)'(1)$  && $(f^*)''(1)$   \\
 \cline{1-6}
  & Normal                                     &&   $-\infty $  &&   -    \\
  & Tukey$(\lambda )$, $\quad 0<\lambda <1$    &&   $-\infty $  &&   -    \\
  & Logistic (Tukey$(0)$)                      &&   $-6$         &&   -    \\
Medium  & Extreme Value                        &&   $-4$         &&   -    \\
   & Laplace                                   &&   $-2$         &&   -    \\
   & Exponential                               &&   $-2$         &&   -    \\
\cline{2-6} &Pareto$(a)$, $\quad 1<a$               &&   $0$     && $+\infty $ \\
       &Lognormal                              &&   $0$     && $+\infty $ \\
 Long  & Tukey$(\lambda )$, $-1 < \lambda <0$  &&   $0$     && $+\infty $ \\
       & Cauchy                                &&   $0$     && $4\pi ^2=39.48 $ \\
       & Tukey$(-1 )$                          &&   $0$     && $33.69   $ \\
       & Pareto$(1)$                           &&   $0$     &&  $6$       \\
\cline{2-6} & Tukey$(\lambda )$, $\quad \lambda <-1$    &&   0       &&  0         \\
Very long      & Pareto$(a)$, $0<a<1$          &&  $0$      &&  0         \\
 \end{tabular}
\end{center}
\end{small}
\end{table}

\subsubsection*{Student's t-distribution with $\nu =2$ degrees of freedom.}

\cite{jones-2002} gives the density quantile function
of Student's t-distribution with $\nu =2$ degrees of freedom,
and by numerical integration one can find the normalizing constant $\kappa = 0.20826$ to obtain its \pdQ \ $f^*(u)=fQ(u)/\kappa =\{2u(1-u)\}^{3/2}/\kappa $. The reader can readily verify that  $f^*(1)=0$, $(f^*)'(1)=0$ and  $(f^*)''(1)=+\infty $ so this $f^*$ has long tails, but shorter than those of the Cauchy distribution. The tail behavior for large $\nu $ may possibly be found using results in \cite{schlut-2012}.

 \subsubsection*{The Pareto distributions.}
The Pareto($a$) family, with $a>0$ has $f^*_a(u)= (2+\frac {1}{a})(1-u)^{1+1/a}\to 0$ as $u\to 1.$
So $(f^*_a)'(u)= -(1+\frac {1}{a})(2+\frac {1}{a})(1-u)^{1/a}\to 0$ as $u\to 1$, again, for all $a>0.$
Therefore $n^*\geq 2$ and the right-hand tails are \lq long\rq\  or \lq very long\rq.\  The second derivative satisfies, as $u\to 1$
\begin{eqnarray*}
    (f^*_a)''(u) &=&  \frac{(1+a)(1+2a)}{a^3}\,(1-u)^{1/a-1}    \\
                   &\to & \left\{
                   \begin{array}{ll}
                     0, & 0 <a < 1 \\
                     6, &   \qquad a=1\\
                   +\infty  , & \qquad 1<a ~.
                   \end{array}
                 \right.
\end{eqnarray*}
Therefore the tails are long ($n^*=2$) if and only if $a\geq 1$ and otherwise they are very long ($n^*\geq 3$). Further, $(f^*_a)^{(n)}(1)=0$ for $0<a<1/2^{n-1}$ for all $n\geq 2$. Thus there exist distributions with tails of every $n^*$-order.

\subsubsection*{The symmetric Tukey($\lambda $) distributions.}
The Tukey$(\lambda )$ family provides a wide range of tail-weights. Starting with $f^*_\lambda $ from
Table~\ref{table1}, the first two derivatives are:
\begin{eqnarray*}
  \kappa _{\lambda }\,f^*_\lambda (u) &=& \{u^{\lambda -1} +(1-u)^{\lambda -1}\}^{-1} \\
   \kappa _{\lambda }\,(f^*_\lambda )'(u) &=& (1-\lambda )(\kappa _{\lambda }\,f^*_\lambda (u) )^2\{u^{\lambda -2}-(1-u)^{\lambda -2}\} \\
   \kappa _{\lambda }\,(f^*_\lambda )''(u)&=& 2(1-\lambda )^2\frac {\{u^{\lambda -2}-(1-u)^{\lambda -2}\}^2}{
   \{u^{\lambda -1}+(1-u)^{\lambda -1}\}^3}+(1-\lambda )(\lambda -2)\frac {\{u^{\lambda -3}+(1-u)^{\lambda -3}\}}{ \{u^{\lambda -1}+(1-u)^{\lambda -1}\}^2} \\
   &\sim &2(1-\lambda )^2\frac {\{1-(1-u)^{\lambda -2}\}^2}{
   \{1+(1-u)^{\lambda -1}\}^3}+(1-\lambda )(\lambda -2)\frac {\{1+(1-u)^{\lambda -3}\}}{ \{1+(1-u)^{\lambda -1}\}^2}~,
\end{eqnarray*}
as $ u\to 1$.   From the first equation above $f^*_\lambda (1)>0$ if and only if $\lambda \geq 1$ and in this case  \lq short\rq\  tails are obtained.   For $\lambda  <1$ examination of the first derivative $(f^*_\lambda )'(u)$ as $u\uparrow 1$ yields the value 0 for   $\lambda <0$,  the value $-6 $ for $\lambda =0$ and $-\infty $ for $0<\lambda < 1;$ in these last two cases $n^*=1$ and the tails are \lq medium\rq.\
  For $\lambda <0$ one can see from the last displayed expression  that $(f^*_\lambda )''(1)=0$ when $\lambda <-1$   and  so the tails are \lq very long\rq.\   For $-1\leq \lambda <0$ the second derivative $(f^*_\lambda )''(1)> 0$ and the
   tails are \lq long\rq.\    In particular, for $\lambda =-1$ one has $(f^*)''(1)\approx 33.7$ which is not far   from  that of the Cauchy $4\pi ^2,$ see Table~\ref{table8}.

\section{Summary and further research}\label{sec:summary}

The \pdQ transformation from $f$ to $f^*$ is quite powerful, and allows for a different look at distributional shapes, on a common finite domain [0,1] where location, scale and gaps are not distractions. This has enabled us to compare discrete with continuous distributions by applying the Hellinger metric and/or Kullback-Leibler divergences to their respective
\pdQ s.  It also facilitated finding \lq closest\rq\  symmetric distributions to a given \pdQ, so that they could be ordered by their distance or divergences from the symmetric class.  Further, we have classified tail behavior using boundary derivatives of the \pdQ s.

With regard to inference, we defined empirical \pdQ s in both the discrete and continuous case. The latter are based
on quantile density estimators of \cite{prst-2016}, and generally require moderately large sample sizes of 500 or more. Given such  a sample, we showed that one could fit a parametric shape model to it by minimizing over the shape parameter the Hellinger distance of the proposed model \pdQ from the empirical \pdQ . This $H$-\pdQ \, method is
location-scale free, and we demonstrated its effectiveness relative to the ppcc method of \cite{fill-1975} for the
symmetric Tukey($\lambda $) family as well as gamma and Weibull models. The new method is a {\em global} density fitting technique, but rather than finding the distance of an empirical density $f_n(x)$ with a proposed $f(x)$ on an infinite domain, which is hard to do, it finds the distance between the normalized $f_n(Q_n(u))$ and a proposed $f(Q(u))$ on [0,1], in effect comparing  the densities $f_n$ and $f$ at their respective quantiles. This procedure
appears more resistant to outliers than traditional fitting methods such as maximum likelihood.

In the asymmetric case one could similarly try fitting the generalized lambda distributions without having to estimate location and scale; background material is in \cite{gilch-2000}, \cite{kardud-2000}, \cite{gld} and \cite{prst-2016}.  Another family with two shape parameters \cite{dagum-1977} might also be
 fitted by the $H$-\pdQ \; technique.

As a specific application, given wool fibre diameter data and hence an empirical \pdQ , we showed that one could globally fit shape families such as Gamma and Weibull to it.  When parametric models do not fit such fibre diameter data well, a non-parametric approach to classifying wool is possible using quantile methods.  Fo example one rcould replace the mean by the median and the coefficient of
variation by an interquantile range divided by the median, which is just the reciprocal of the
standardized median investigated by \cite{S-2013}.

With regard to theoretical research on \pdQ s, what happens with another
application of the \pdQ -transformation to $f^{**}=(f^*)^*$?  And further, with $n$ iterations $f^{(n+1)*}=(f^{n*})^*$ for $n\geq 2$, and $n\to \infty $ ? Preliminary work suggests that under weak
conditions, such as the $n$th power integrability of $f$ for all $n$,
 the limit $\lim  _{n\to \infty}\,f^{(n)*}(u)$ exists and is the \lq shapeless\rq\  uniform distribution on [0,1].

While the \pdQ s of discrete distributions are of interest, we have not delved into them much here.
For example, the moments in Table~\ref{table2} could  include those of discrete distributions.
 It would also be of interest to extend the tail-weight analysis of Section~\ref{sec:tailwt} to discrete distributions, although the boundary derivatives
could not be found as limits of derivatives near the boundary.

With regard to right tail-weight, those \pdQ s of $n^*$-order, but infinite $n^*$ derivative at 1, could be compared by looking at the relative rates at which these $n^*$ derivatives approach infinity. It would be of interest to include  Student's $t_\nu $ \pdQ s for $\nu >1$ degrees of freedom, as well as  many other long and very long tailed distributions, in Table~\ref{table8}.

Extending the concept of \pdQ s to the  multivariate case appears feasible, if challenging, and likely to lead to numerous applications for multivariate data, because
the metrics, divergences and empirical \pdQ s are still available.

\bigskip

{\em : Acknowledgments:  The author is indebted to the Editors and the referees whose detailed comments
and suggestions greatly improved the clarity of this manuscript.}

\clearpage
\newpage

\section{Appendix}\label{sec:appendix}

\subsection{Wool fibre diameter data}\label{sec:appenda}
\begin{figure}[t!]
\begin{center}
\includegraphics[scale=.7]{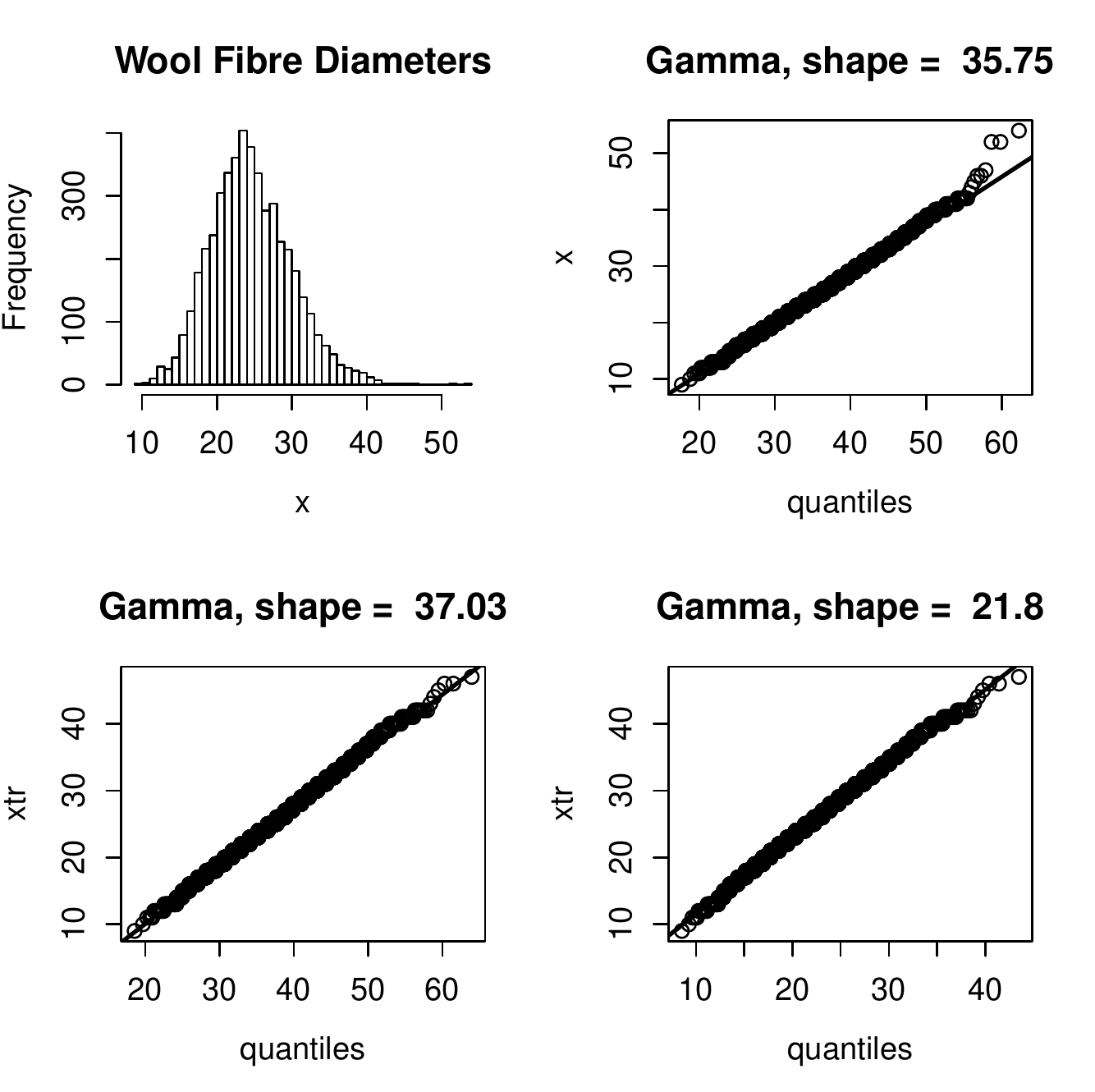}
\caption{\footnotesize  \em Histogram of wool fibre diameter data $x$ and QQ-plots of fitted gamma model
for these data and also for two fitted gamma models for truncated data $xtr$\,; data and details of analysis are found in Section~\ref{sec:appenda}. \label{fig6}}
\end{center}
\end{figure}
A histogram of wool fibre diameters is shown in the upper left hand plot of Figure~\ref{fig6};
it is reproduced from \cite[p.246]{botha-2010}. The data can be loaded into R and a statistical
summary obtained with the following commands: \  here $d=9,10,\dots ,54$ is a vector of diameters in microns, $f$ is a vector
of their frequencies and $x$ is the sorted data, numbering $n=4817.$
\begin{small}
\begin{verbatim}
d <- c(seq(9,54))
f <- c(1,1,3,10,29,25,43,79,117,178,216,238,305,337,361,404,378,336,277,288,227)
f <- c(f,215,181,139,113,79,62,48,31,27,22,19,12,7,1,1,1,2,1,0,0,0,0,2,0,1)
x <- rep.int(d,f)
summary(x)
  Min.    1st Qu.     Median     Mean     3rd Qu.    Max.
  9.00     21.00      25.00      25.08     28.00    54.00
\end{verbatim}
\end{small}
The top right plot of Figure~\ref{fig6} shows a QQ plot of these data versus the quantiles of a Gamma
distribution with shape parameter $35.75$, the best fitting model of Table~\ref{table6}. Visually it
is very similar to QQ-plots using the Gamma model and shape estimates with values ranging from 22 to 37.
The QQ-plots for the best-fitting Weibull models in Table~\ref{table6} are not shown because they
are non-linear and so the Weibull models are  unsuitable for fitting these data.

It was thought that the outliers 52, 52 and 54 might be affecting the fits of the Gamma models, so
they were fitted again after omitting them; the truncated data is called \lq xtr\rq,\ and it has
maximum 47, one sd less than the omitted values. In the bottom left-hand plot of Figure~\ref{fig6}
it is seen that the QQ-plot for the best fitting Gamma model found by method $H$-\pdQ , is almost
linear, and similarly on its right for shape obtained by the MLE method.
\clearpage
\newpage
\subsection{Proof of Proposition~\ref{prop1}}\label{sec:appendb}
\begin{figure}[t!]
\begin{center}
\includegraphics[scale=.7]{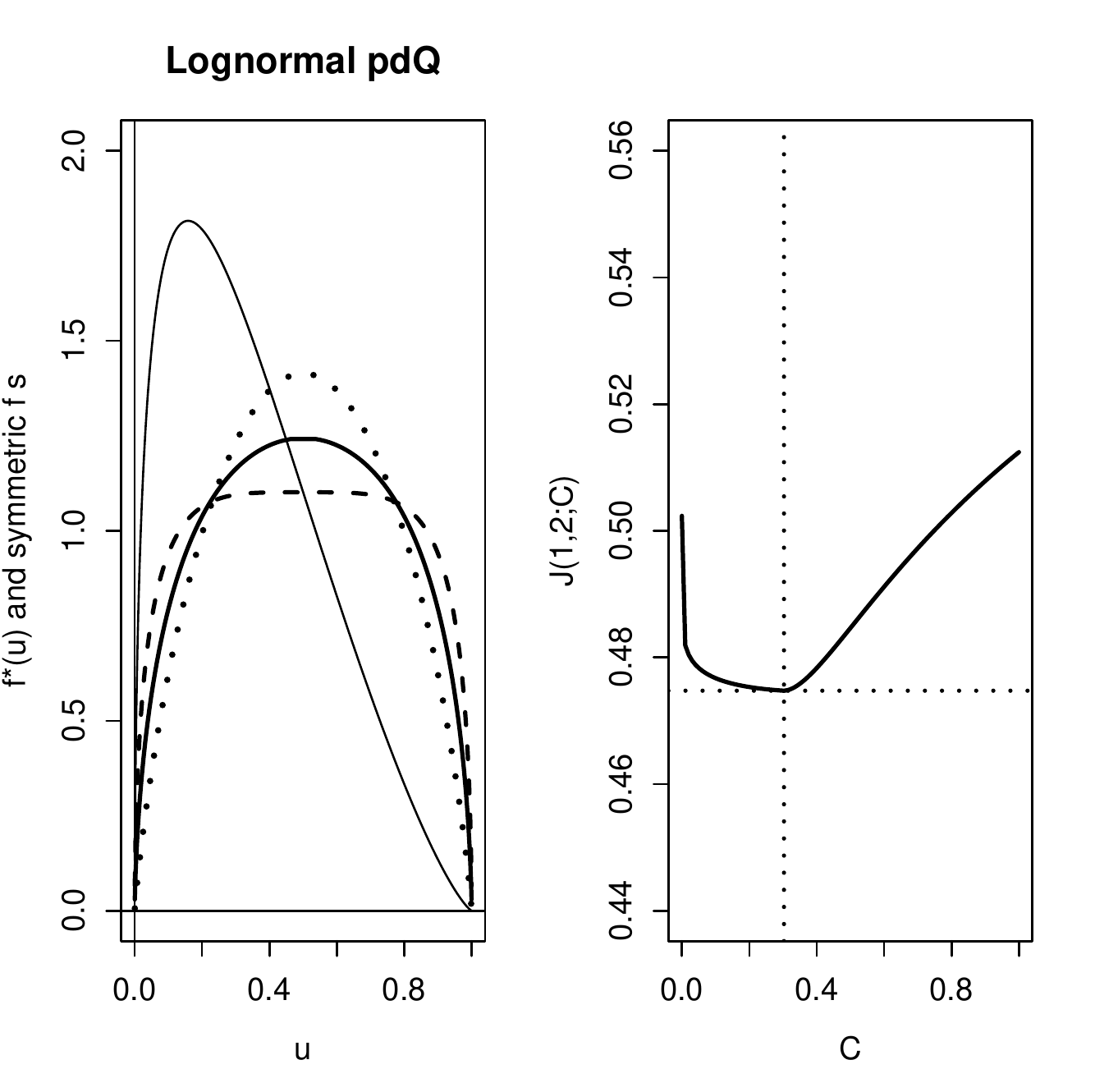}
\caption{\footnotesize  \em  On the left is shown the lognormal \pdQ  $f^*(u)$ (thin solid line) and the respective closest symmetric densities on [0,1] that minimize $I(1:2)$ (dashed line), $I(2:1)$ (dotted line) and $J(1,2)$
(thick solid line). On the right is shown the graph of $J(1,2; C)$ for various $C$; the $C_\text{opt}=0.303$ which minimizes the \kld distance $J(1,2)$ is marked by a vertical line. The closest symmetric density in the Hellinger metric to the lognormal \pdQ is the same as that minimizing $J(1,2)$. \label{fig7}}
\includegraphics[scale=.7]{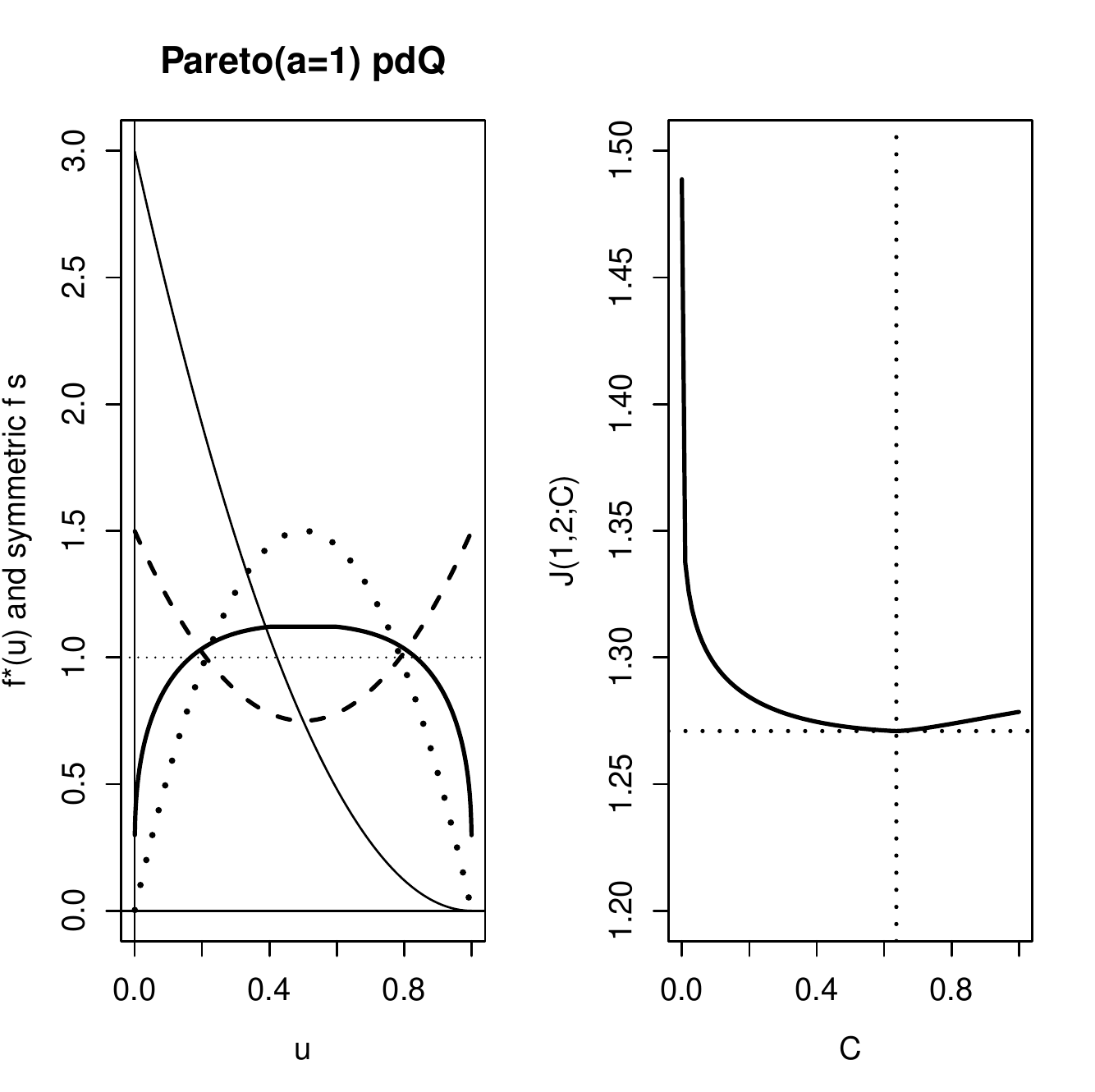}
\caption{\footnotesize On the left is shown the Pareto($a=1$) \pdQ  $f^*(u)$ (thin solid line) and the respective closest symmetric densities on [0,1]  that minimize $I(1:2)$ (dashed line), $I(2:1)$ (dotted line) and $J(1,2)$
(thick solid line). On the right is shown the graph of $J(1,2; C)$ for various $C$; the $C_\text{opt}=0.636$ which minimizes the \kld distance $J(1,2)$ is marked by a vertical line.
The dotted horizontal line at 1 marks the uniform density, which is the closest symmetric density on [0,1] in the Hellinger metric to this Pareto \pdQ ,\ and for this example it differs from the symmetric density minimizing $J(1,2)$. \label{fig8}}
\end{center}
\end{figure}

{\em Part (a) is Theorem 2.2 of  \cite{with-2010}.
Part (b) is a straightforward modification of the proof of Part(a):\  let
$X$ have a distribution on the integers with $P(X=i)=p_i$ and let $Y$ be a symmetric distribution on the integers with
$P(Y=i)=q_i$, so $q_{-i}=q_i$, $i\geq 1.$  Amongst such $q$ we will minimize the  divergence $I(p:q)=\sum _i p_i\,\ln (p_i/q_i)=\sum _i p_i\ln (p_i)-p_0\ln(q_0)-\sum _{i=1}^{\infty }(p_i+p_{-i})\,\ln (q_i)$. Taking $\lambda $ as a Lagrange multiplier, we need to solve $0\equiv \frac {\partial}{\partial q_i} \{I(p:q)+\lambda \sum _iq_i \}$. This yields $0=-p_0/q_0+\lambda $ and $0=-(p_i+p_{-i})/q_i+2\lambda $ for $i>0$. Therefore $q_0=p_0/\lambda $ and $q_i=(p_i+p_{-i})/(2\lambda )$. The condition $\sum _iq_i=1$ then implies $\lambda =1$.  Thus $q_i=(p_i+p_{-i})/2$ for all $i$.
Clearly this result can be extended to an arbitrary discrete distribution on a lattice of points ${\cal X}=\{x_i\}$ by extending ${\cal X}$ to ${\cal X}\cup -{\cal X}$, letting $ P(X=x_i)=p_i$ and using the same argument. The absolutely
continuous case then follows by approximating the density over a lattice with smaller and smaller increments and taking the limit.

For part (c), again return to a given discrete distribution $p$ on the integers, and a symmetric distribution $q$ on
the integers. For each positive integer $i$ introduce $2a_i=\ln(p_ip_{-i})$ and $\nu _i =\sqrt {p_ip_{-i}}\,=\exp(a_i).$  Further let $\overline p_i=(p_i+p_{-i})/2.$ To minimize the \kld between
$p$ and $q$ by choice of symmetric $q$ we need to minimize
\[J(p,q) =\sum _i p_i\ln (p_i)-p_0\ln(q_0)-2\sum _{i=1}^{\infty }\overline p_i\,\ln (q_i)
+ q_0\ln (q_0/p_0)+2 \sum _{i=1}^\infty q_i(\ln (q_i)-a_i)~.\]
Setting the derivatives $0\equiv \frac {\partial}{\partial q_i} \{J(p,q)+\lambda \sum _iq_i \}$ yields
$0=-p_0/q_0+\ln(q_0/p_0)+1+\lambda $ and $0= 2\{-\overline p_i/q_i +\ln(q_i)-a_i+1+\lambda \}$ for $i>0$.
The solution for $q$ is implicit in:
\begin{eqnarray}\label{klditer}\nonumber
    q_0 &=& p_0\exp \{p_0/q_0\}\;e^{-\lambda -1}~;       \\
    q_i &=& \nu _i\,\exp \{\overline p_i/q_i\}\;e^{-\lambda -1} ~, \text { for } i\geq 1~;\\ \nonumber
      1 &=& q_0 +2 \sum _{i=1}^\infty q_i       ~.
\end{eqnarray}
As in parts (a) and (b), the proof of part (c) is completed by a limiting argument.

\smallskip
{\bf Algorithm for computing the solution to (\ref{klditer}).}
Given a discrete distribution $p$, ${\cal X}$, one may solve the equations (\ref{klditer}) iteratively. For the
continuous case we solved the equations (\ref{Jmin}) as follows:
\begin{enumerate}
  \item  Given $f_1$ compute $\nu (x)= \sqrt {f_1(x)\,f_1(-x)}\,$ and $\overline f_1(x)\equiv \{f_1(x)+f_1(-x)\}/2$.
  \item  Fix $C=e^{-\lambda -1}>0$ and for each $x$ in a fine grid over the support of $f_1$ solve for $\beta (x;C)$
\[\beta (x) = \nu (x) \,\exp \left \{\overline f_1(x)/\beta (x)\right \}\,C ~.\]
 Numerically compute $d_C=\int \beta (x;C)\,dx $, normalize $\beta (x;C)$ to a probability density $f_2(x;C)=\beta (x;C)/d_C$ and use it to calculate $J(1,2; C)$.
\item   Repeat the last step for a range of $C$ values to locate the $C=C_\text{opt}$ for which $J(1,2;C)$ is minimized. This $f_2(x;C_\text{opt})$ is the solution guaranteed by (\ref{Jmin}).
\end{enumerate}
}
{\bf Examples of computing (\ref{klditer}).}
In the left plot of Figure~\ref{fig7} are shown the symmetric densities closest to the lognormal \pdQ \;(in minimizing the divergences $I(1:2)$, $I(2:1)$ and $J(1,2)$). Also shown are values $J(1,2;C) = \int  \{f_1(u)\log\left (f_1(u)/f_C(u)\right)+f_C(u)\log\left(f_C(u)/f_1(u)\right )\}\, du ,$ where $f_C(u)\equiv f_2(x;C)=\beta (x;C)/d_C$ in the iterative Step 2 above. Figure~\ref{fig8} gives the
corresponding results for the Pareto $a=1$. An R script for
creating Figure~\ref{fig8} is in the supplementary online material.

\end{document}